\documentclass[12pt]{article}
\usepackage{amsmath}
\usepackage{amssymb}
\usepackage[cp1251]{inputenc}
\usepackage[russian]{babel}
\newcommand{\il}[2]{\int\limits_{#1}^{#2}}

\newcommand{\ph}{\phantom{a}}
\newcommand{\phh}{\phantom{aaa}}

\newcommand{\sist}[2]{\left\{
\begin{array}{l}
{#1}\\
\ph\\
{#2}
\end{array}
\right.}
\newcommand{\sisttt}[3]{\left\{
\begin{array}{l}
{#1}\\
\ph\\
{#2}\\
\ph\\
{#3}
\end{array}
\right.}

\begin{document}

MSC 34D20

\vskip 20pt

\centerline{\bf  Comparison criterion for second order}
\centerline{\bf  Riccati equations}

\vskip 10 pt

\centerline{\bf G. A. Grigorian}

\vskip 10 pt

\centerline{0019 Armenia c. Yerevan, str. M. Bagramian 24/5}
\centerline{Institute of Mathematics NAS of Armenia}
\centerline{E - mail: mathphys2@instmath.sci.am, \ph phone: 098 62 03 05, \ph 010 35 48 61}

\vskip 20 pt
\noindent
Abstract.
Three comparison criteria are obtained for second order Riccati equations. On the basis of these criteria some global existence theorems are proved mentioned equations. The results obtained are used to derive a non oscillation criterion for three dimensional linear systems of ordinary differential equations.


\vskip 20 pt

{\bf 1. Introduction}. Let $a(t), \ph b(t), \ph c(t), \ph d(t)$ and $e(t)$ be real-valued continuous functions on $[t_0,+\infty)$. Consider the second order Riccati equation
$$
y'' + 3 a(t) y y' + b(t) y' + a^2(t) y^3 + c(t) y^2 + d(t) y + e(t) = 0, \phh t \ge t_0. \eqno (1.1)
$$
Throughout we will assume that $a(t), \ph b(t) \in C^1([t_0,\infty)).$ This equation serves as an important tool for the study of properties of solutions of first order three dimensional equations, in particular, for the study of properties of solutions of third order linear ordinary differential equations. One of ways to start to apply this equation to the mentioned above systems is to obtain comparison criteria for it. It should be noticed here that  this approach has been used to first order Riccati equations (see [1,2]) allowing to obtain several results for qualitative study of solutions of some types of equations (see e. g. [3--13]).

In this paper we prove two comparison criteria for second order Riccati equations. We use these criteria to obtain some global solvability theorems for Eq. (1.1) and some nonoscillation criteria for first order three dimensional linear systems of ordinary \linebreak differential equations.

\vskip 10pt

{\bf 2. Auxiliary propositions.} Let $a_1(t), \ph b_1(t), \ph c_1(t), \ph d_1(t)$ and $e_1(t)$ be real-valued continuous functions on $[t_0,+\infty)$. Along with (1.1) consider the equation
$$
y'' + 3 a_1(t) y y' + b_1(t) y' + a_1^2(t) y^3 + c_1(t) y^2 + d_1(t) y + e_1(t) = 0, \phh t \ge t_0. \eqno (2.1)
$$
We set $\nu(t,u,v,u_1,v_1) \equiv u_1 - v_1 +\frac{3}{2} a(t) (u^2 - v^2) + b(t) (u - v), \ph \Gamma(t,u,v) \equiv a^2(t)[u^2 + uv + v^2] + (c(t) - \frac{3}{2}a'(t))(u + v) - b'(t) + d(t), \ph J(t,u,v)\equiv (u - v)\Gamma(t,u,v), \ph L(t,u,v) \equiv 3(a_1(t) - a(t))) uv  + (b_1(t) - b(t)) v + (a_1^2(t) - a^2(t)) u^3 + (c_1(t) - c(t)) u^2 + (d_1(t) - d(t)) u + e_1(t) - e(t), \ph t \ge t_0, \ph u, v, u_1, v_1 \in \mathbb{R.}$ Let $y_0(t)$ and $y_1(t)$ be solutions of Eq. (1.1) and Eq. (2.1) respectively on $[t_1,t_2) \subset [t_0,\infty)$. Then it is not difficult to  verify that

\noindent
$[y_0(t) - y_1(t)]'' + \frac{3}{2} a(t) [(y_0(t) - y_1(t))(y_0(t) + y_1(t))]' + b(t)[y_0(t) - y_1(t)]' + a^2(t) [y_0(t) - y_1(t)] [y_0^2(t) + y_0(t) y_1(t) + y_1^2(t)] + c(t)(y_0(t) - y_1(t))(y_0(t) + y_1(t)) + d(t) (y_0(t) - y_1(t)) - L(t,y_1(t),y_1'(t)) = 0, \ph t \in [t_1,t_2).$ Let us integrate this equality from $t_1$ to $t$. After using the rule of integration by parts and making some simplifications we obtain
$$
[y_0(t) - y_1(t)]' + \Bigl\{\frac{3}{2} a(t) (y_0(t) + y_1(t)) + b(t)\Bigr\}(y_0(t) - y_1(t)) - \nu(t, y_0(t),y_1(t), y_0'(t),y_1'(t)) -
$$
$$
- \il{t_1}{t}J(\tau,y_0(\tau),y_1(\tau)) d \tau   -\il{t_1}{t} L(\tau, y_1(\tau),y_1'(\tau))d \tau = 0, \ph t \in [t_1,t_2). \eqno (2.2)
$$
It is clear from here that $y_0(t) - y_1(t)$ is a solution of the linear equation
$$
x' + \Bigl\{\frac{3}{2} a(t) (y_0(t) + y_1(t)) + b(t)\Bigr\}x  = \nu(t, y_0(t),y_1(t), y_0'(t),y_1'(t)) +
$$
$$
 +\il{t_1}{t}J(\tau,y_0(\tau),y_1(\tau)) d \tau   + \il{t_1}{t} L(\tau, y_1(\tau),y_1'(\tau))d \tau = 0, \ph t \in [t_1,t_2).
$$
Then by the Cauchy formula we have
$$
y_0(t) - y_1(t) = \exp\biggl\{-\il{t_1}{t}\biggl\{\frac{3}{2} a(\tau)(y_0(\tau) + y_1(\tau)) + b(\tau)\biggr\} d \tau\biggr\}\times
$$
$$
\times \Biggl[y_0(t_1) - y_1(t_1) + \il{t_1}{t}\exp\biggl\{\il{t_1}{\tau}\biggl[\frac{3}{2}a(s)(y_0(s) + y_1(s)) + b(s)\biggr] d s\biggr\} \times
$$
$$
\times\biggl(\nu(t_1,y_0(t_1),y_1(t_1),y_0'(t_1),y'_1(t_1)) (\tau - t_1) + \il{t_1}{\tau}[J(s,y_0(s),y_1(s)) +
$$
$$
+L(s,y_1(s),y_1'(s))]d s\biggr) d \tau\biggr], \phh t \in [t_1,t_2).  \eqno (2.3)
$$
We set $D(t) \equiv 2(2 c(t) - 3 a'(t))^2 + a^2(t)(d(t) - b'(t)). \ph t \ge t_0$.

\vskip 10pt

{\bf Lemma 2.1.} {\it If $a(t) \ne 0, \ph D(t) \ge 0, \ph t \ge t_0$, then $\Gamma(t,u,v) \ge 0, \ph t \ge t_0, \ph u, v \in \mathbb{R}.$
}

Proof. We have
$$
\sist{\frac{\partial \Gamma(t,u,v)}{\partial u} = a^2(t)(2 u + v) + c(t) - \frac{3}{2} a'(t),}{\frac{\partial \Gamma(t,u,v)}{\partial u} = a^2(t)(u + 2 v) + c(t) - \frac{3}{2} a'(t), \ph t \in[t_1,t_2), \ph u,v \in \mathbb{R}.}  \eqno (2.4)
$$
Since $a(t) \ne 0, \ph t \in [t_1,t_2)$ one can easily verify that the functions $u_0(t)  = v_0(t) = \frac{3 a'(t) 2 c(t)}{2 a^2(t)}, \ph t \in [t_1,t_2)$ form the unique solution $(u_0(t),v_0(t))$ of the system
$$
\sist{ a^2(t)(2 u + v) = \frac{3}{2} a'(t)  - c(t),} {a^2(t)(u + 2 v)  = \frac{3}{2} a'(t) - c(t), \ph t \in[t_1,t_2),  \ph u,v \in \mathbb{R}.}
$$
Then by virtue of (2.4) it follows from the conditions of the lemma that
$$
\min\limits_{u,v\in\mathbb{R}}\Gamma(t,u,v) = \Gamma(t,u_0(t),v_0(t)) = \frac{D(t)}{a^2(t)} \ge 0, \ph t \in [t_1,t_2).
$$
Therefore, $\Gamma(t,u,v) =  \frac{D(t)}{a^2(t)} \ge 0, \ph t \in [t_1,t_2).$ The lemma is proved.

Let $X(t), \ph Y(t), \ph Z(t)$ and $W(t)$ be real-valued continuous functions on $[t_0,\infty)$ and let $X(t)\ne 0, \ph t \ge t_0, \ph X(t) \in C^1([t_0,\infty)).$ Consider the linear system
$$
\sisttt{\phi' = \phantom{aaaaaaa} a(t) \psi,}{\psi' = \phantom{aaaaaaaaaaaaaa} X(t) \chi,}{\chi' = Y(t)\phi + Z(t)\psi + W(t)\chi, \ph t \ge t_0.} \eqno (2.5)
$$
The substitution
$$
\psi = y \phi \eqno (2.6)
$$
reduces this system into
$$
\sisttt{\phi' = \ph a(t) y \phi,}{[y' + a(t) y^2]\phi =  \ph X(t) \chi,}{\chi' = [Y(t)\phi + Z(t) y]\phi + W(t)\chi, \ph t \ge t_0.} \eqno (2.7)
$$
Since $X(t) \ne 0, \ph t \ge t_0,$ from here we get
$$
\chi = \frac{y' + a(t) y^2}{X(t)} \phi, \phh t \ge t_0. \eqno (2.8)
$$
This together with the last equation of the system (2.7) implies
$$
\Bigl(\frac{1}{X(t)}\Bigr)'(y' + a(t) y^2)\phi + \frac{1}{X(t)}[y'' + a'(t) y' + 2 a(t) y y']\phi = [Y(t) + Z(t) + W(t) \frac{y' + a(t) y^2}{X(t)}] \phi,
$$
$t \ge t_0$. After some simplifications from here we derive
$$
y'' + 3 a(t) y y' - \Bigl[\frac{X'(t)}{X(t)} + W(t)\Bigr] y' + a^2(t) y^3  + \Bigl[W(t) a(t) - \frac{X'(t)}{X(t)} a(t) + a'(t)\Bigr] y^2 + X(t) Z(t) y +
$$
$$
+ X(t) Y(t) = 0, \phh t \ge t_0. \eqno (2.9)
$$
This equation coincides with Eq. (1.1) provided
$$
\sisttt{\frac{X'(t)}{X(t)} + W(t) = - b(t),}{W(t) a(t) - \frac{X'(t)}{X(t)} a(t) + a'(t) = c(t),}{X(t) Z(t) = d(t), \ph X(t) Y(t) = e(t), \ph t \ge t_0.}
$$
Assume $a(t) \ne 0, \ph t \ge t_0$. Then the last system is equivalent to the following one.
$$
\sisttt{X'(t) = \Bigl[\frac{a'(t)}{a(t)} - \frac{c(t)}{2 a(t)} - \frac{b(t)}{2}\Bigr] X(t),}{W(t) = -\frac{X'(t)}{X(t)} - b(t),}{Y(t) = \frac{e(t)}{X(t)}, \ph Z(t) = \frac{d(t)}{X(t)}, \ph t \ge t_0.}
$$
Hence, if in particular, $X(t) = \exp\biggl\{\il{t_0}{t}\Bigl[\frac{a'(\tau)}{2 a(\tau)} - \frac{c(\tau)}{2 a(\tau)} - \frac{b(\tau)}{2}\Bigr]d \tau\biggr\}, \ph Y(t) = \frac{e(t)}{X(t)}, \ph Z(t) = \frac{d(t)}{X(t)}, \ph W(t) = \frac{b(t)}{2} + \frac{c(t)}{2 a(t)} - \frac{a'(t)}{a(t)}, \ph t \ge t_0,$ then Eq. (2.9) coincides with Eq. (1.1). It follows from here and from (2.6) - (2.8) that all solutions $y(t)$ of Eq. (1.1), existing on any interval $[t_1,t_2) \subset [t_0,\infty)$, are connected with solutions $(\phi(t),\psi(t),\chi(t))$ of the system (2.5) by relations
$$
\sist{\phi(t) = \phi(t_1) \exp\biggl\{\il{t_1}{t} a(\tau) y(\tau) d \tau\biggr\}, \phi(t_1) \ne 0,}{\psi(t) = y(t) \phi(t), \chi(t) = \frac{y'(t) + a(t) y^2(t)}{X(t)}, \ph t \in [t_1,t_2).} \eqno (2.10)
$$

{\bf Definition 2.1.} {\it An interval $[t_1,t_2) \subset [t_0,\infty)$ is called the maximum existence interval for a solution $y(t)$ of Eq. (1.1), if $y(t)$ exists on $[t_1,t_2)$ and cannot be continued to the right from $t_2$ as a solution of Eq. (1.1).}

{\bf Lemma 2.2.} {\it Let $a(t) \ne 0, \ph t \in [t_1,t_2) \subset[t_0,\infty)$ and let $y(t)$ be a solution of Eq. (1.1) on $[t_1,t_2)$. If the function $F(t) \equiv \il{t_0}{t} a(\tau) y(\tau) d \tau, \ph t \in [t_1,t_2)$ is bounded from below on $[t_1,t_2)$, then $[t_1,t_2)$ is not the maximum existence interval for $y(t)$.
}

Proof. Since $a(t) \ne 0, \ph t \in [t_1,t_2)$ by (2.10) the functions $\phi(t) \equiv \exp\biggl\{\il{t_1}{t} a(\tau) y(\tau) d \tau\biggr\}, \linebreak  \psi(t) \equiv y(t)\phi(t), \ph \chi(t) \equiv \frac{y'(t) + a(t) y^2(t)}{X(t)}, \ph t \in [t_1,t_2)$  form a solution $(\phi(t),\psi(t),\chi(t))$ of the system (2,5) on $[t_1,t_2)$, which is continuable on $[t_1,\infty)$ as a solution of the system (2.5) on $[t_1,\infty)$. It follows from here and from the boundedness from below of $F(t)$ that $\phi(t) \ne 0, \ph t \in [t_1,t_3)$ for some $t_3 > t_2.$ Then by (2.6) $\widetilde{y}(t) \equiv \frac{\psi(t)}{\phi(t)}, \ph t \in [t_1,t_3)$ is a solution of Eq. (1.1) on $[t_1,t_3)$. Obviously $\widetilde{y}(t) = y(t), \ph t \in [t_1,t_2)$. Hence, $\widetilde{y}'(t) = y'(t), \ph t \in [t_1,t_2)$. It follows from here that $[t_1,t_2)$ is not the maximum existence interval for $y(t)$. The lemma is proved.

Consider the differential inequality
$$
\eta' + 3 a(t) \eta \eta' + b(t) \eta' + a^2(t)\eta^3 + c(t)\eta^2 + d(t) \eta + e(t) \ge 0, \ph t\ge t_0. \eqno (2.11)
$$

{\bf Remark 2.1.} {\it It is not difficult to verify that to verify that if $a(t)\ne 0,  \ph t \in [t_0,T]$ for some $T > t_0$, then $\eta_\lambda \equiv \lambda + \max\limits_{t\in[t_0,T]}\frac{|c(t)| + |d(t)| + |e(t)|}{a^2(t)}, \ph t \in [t_0,T],$ where $\lambda = const \ge 0$, is a solution of the inequality (2.11), and if $a(t) \ne 0, \ph t \ge t_0, \ph M\equiv \sup\limits_{t\ge t_0} \frac{c(t)| + |d(t)| + |e(t)|}{a^2(t)} < \infty$, them $\eta_\lambda(t) \equiv \lambda + M, \ph t \ge t_0$, is a solution of the inequality (2.11) on $[t_0,\infty)$.}

{\bf Lemma 2.3.} {\it Let $y(t)$ be a solution of Eq. (1.1)  on $[t_1,t_2)$ and let $\eta(t)$ be a solution of the inequality (2.11) on $[t_1,t_2)$ such that $\eta(t_1) \ge y(t_1), \ph \nu(t_1,\eta(t_1),y(t_1)\eta'(t_1),y(t_1))\ge 0.$ If $a(t) \ne 0, \ph D(t) \ge 0, \ph t \in [t_1,t_2)$, then
$$
\eta(t) \ge y(t), \ph t \in [t_1,t_2), \eqno (2.12)
$$
$$
\nu(t,\eta(t),y(t),\eta'(t),y'(t))\ge 0, \ph t\in [t_1,t_2), \eqno (2.13)
$$
}

Proof. We have
$$
[\eta(t) - y(t)]'' +\frac{3}{2} a(t)[(\eta(t) + y(t))(\eta(t) - y(t))]' + b(t)(\eta(t) - y(t))' +
 $$
 $$
 +a^2(t)(\eta(t) - y(t))[\eta^2(t) + \eta(t) y(t) + y^2(t)] + c(t)(\eta(t) + y(t))(\eta(t) - y(t)) +
$$
$$
+  d(t)(\eta(t) - y(t)), \ph t \in [t_1,t_2).
$$
Integrating this inequality from $t_1$ to $t$ and making some simplifications we obtain
$$
\mu(t)\equiv [\eta(t) - y(t)]' + \Bigl\{\frac{3}{2} a(t)(\eta(t) + y(t)) + b(t)\Bigr\}(\eta(t) - y(t)) -
$$
$$
-\nu(t_1,\eta(t_1),y(t_1)\eta'(t_1),y(t_1)) - \il{t_1}{t} J(\tau,\eta(\tau),y(\tau)) d\tau \ge 0, \ph t \in [t_1,t_2). \eqno (2.14)
$$
It is clear from here that $\eta(t) - y(t)$ is a solution of the linear equation
$$
x' + \frac{3}{2} a(t)[(\eta(t) + y(t)) x - \nu(t_1,\eta(t_1),y(t_1)\eta'(t_1),y(t_1)) - \il{t_1}{t} J(\tau,\eta(\tau),y(\tau)) d\tau - \mu(t),
$$
$t \in [t_1,t_2).$ Then by the Cauchy formula we have
$$
\eta(t) - y(t)= \exp\biggl\{-\il{t_1}{t}\Bigl[\frac{3}{2} a(\tau)(\eta(\tau) + y(\tau)) + b(\tau)\Bigr]d\tau\biggr\}\times
$$
$$
\times\biggl[\eta(t_1) - y(t_1) + \il{t_1}{t}\exp\biggl\{\il{t_1}{\tau}\Bigl[\frac{3}{2}a(s)(\eta(s) + y(s)) + b(s)\Bigr]ds\biggr\} \times
$$
$$
\times\biggl((\tau - t)\nu(t_1,\eta(t_1),y(t_1),\eta'(t_1),y'(t_1)) + \il{t_1}{\tau}[J(s,\eta(s),y(s)) + \mu(s)]ds\biggr)d\tau\biggr],  \eqno (2.15)
$$
$ t \in [t_1,t_2).$ By Lemma 2.1 it follows from the conditions $a(t) \ne 0, \ph D(t) \ge 0, \linebreak t \in [t_1,t_2)$ that $J(t,\eta(t),y(t)) \ge 0, \ph t \in [t_1,t_2).$ Moreover, it follows from (2.14) that $\mu(t) \ge~ 0, \ph t \in[t_1,t_2)$. Hence, under the initial conditions of the lemma it follows from (2.15) the inequality (2.12). The inequality (2.13) follows immediately  from (2.14), from the initial condition $\nu(t_1,\eta(t_1),y(t_1),\eta'(t_1),y'(t_1)) \ge 0$ and from the inequality  $J(t,\eta(t),y(t)) \ge~ 0, \ph t \in [t_1,t_2).$ The lemma is proved.

Consider the differential inequality
$$
\zeta' + 3 a(t) \zeta \zeta' + b(t) \zeta' + a^2(t)\zeta^3 + c(t)\zeta^2 + d(t) \zeta + e(t) \le 0, \ph t\ge t_0. \eqno (2.16)
$$

{\bf Remark 2.2.} {\it It is not difficult to verify that to verify that if $a(t)\ne 0,  \ph t \in [t_0,T]$ for some $T > t_0$, then $\zeta_\lambda \equiv -\lambda - \max\limits_{t\in[t_0,T]}\frac{|c(t)| + |d(t)| + |e(t)|}{a^2(t)}, \ph t \in [t_0,T],$ where $\lambda = const \ge 0$, is a solution of the inequality (2.16), and if $a(t) \ne 0, \ph t \ge t_0, \ph M\equiv \sup\limits_{t\ge t_0} \frac{c(t)| + |d(t)| + |e(t)|}{a^2(t)} < \infty$, then $\zeta_\lambda(t) \equiv -\lambda - M, \ph t \ge t_0$, is a solution of the inequality (2.16) on $[t_0,\infty)$.}

By analogy with the proof of Lemma 2.3 one can prove the following lemma.

{\bf Lemma 2.4.} {\it Let $y(t)$ be a solution of Eq. (1.1) on $[t_1,t_2)$ and let $\zeta(t)$ be a solution of the inequality (2.16) on $[t_1,t_2)$ such that
$$
\zeta(t_1) \le y(t_1), \ph \nu(t_1,\zeta(t_1),y(t_1),\zeta'(t_1),y'(t_1)) \le 0.
$$
If $a(t) \ne 0, \ph D(t) \ge 0, \ph t \in [t_1,t_2)$, then
$$
\zeta(t) \ge y(t),  \phh
\nu(t,\zeta(t),y(t),\zeta'(t),y'(t))\le 0, \phh t\in [t_1,t_2).
$$
}

\phantom{aaaaaaaaaaaaaaaaaaaaaaaaaaaaaaaaaaaaaaaaaaaaaaaaaaaaaaaaaaaaaaaaaaaaa}$\blacksquare$

Consider the nonlinear system
$$
Y' = F(t,Y), \phh t \ge t_0. \eqno (2.17)
$$
Every solution $Y(t) = Y(t,t_0,Y_0)$ of this system exists either only a finite interval $[t_0,T)$ or is continuable on $[t_0,\infty)$

\vskip 10pt

{\bf Lemma 2.5([14, p. 204, Lemma]).} {\it If a solution $Y(t)$ of the system (2.17) exists only on a finite interval $[t_0,T)$, then
$$
||Y(t)|| \to \infty \ph \mbox{as} \ph t \to T- 0,
$$
where $||Y(t)||$ is any euclidian norm of $Y(t)$ for every fixed $t \in [t_0,T)$.
}

\phantom{aaaaaaaaaaaaaaaaaaaaaaaaaaaaaaaaaaaaaaaaaaaaaaaaaaaaaaaaaaaaaaaaaaaaa}$\blacksquare$

\vskip 10pt

{\bf 3. Comparison criteria.} In this section we  use the results of the previous section to derive three comparison criteria for second order Riccati equations. These criteria we use in the next section to obtain some global solvability criteria for Eq. (1.1).

{\bf Theorem 3.1.} {\it Let the following conditions be satisfied

\noindent
(1) $a(t) > 0, \ph t \in [t_1,t_2),$

\noindent
(2) $D(t) \ge 0, \ph t \in [t_1,t_2),$

\noindent
(3) Eq. (2.1) has a solution $y_1(t)$ on $[t_1,t_2)$ such that
$$
\gamma - y_1(t_1) + \il{t_1}{t}\exp\biggl\{\il{t_1}{\tau}\Bigl[\frac{3}{2} a(s)(\eta(s) + y_1(s)) + b(s)\Bigr]ds\biggr\}\biggl(\il{t_1}{\tau}L(s,y_1(s),y'_1(s)) d s\biggl) d\tau \ge 0,
$$
$t \in[t_1,t_2)$ for some $\gamma \ge y(t_1)$, where $\eta(t)$ is a solution of the inequality (2.11) on $[t_1,t_2)$ with $\eta(t_1) \ge \gamma$.

\noindent
Then every solution $y(t)$ of Eq. (1.1) with the initial conditions $y(t_1) \in [\gamma,\eta(t_1)], \linebreak \nu(t_1,\eta(t_1),y(t_1),\eta'(t_1),y'(t_1)) >~0, \ph \nu(t_1,y(t_1),y_1(t_1),y'(t_1),y_1'(t_1)) \ge 0$, exists on $[t_1,t_2)$ and
$$
y_1(t) \le y(t) \le \eta(t), \ph t \in [t_1,t_2), \eqno (3.1)
$$
$$
\nu(t,\eta(t),y(t),\eta'(t),y'(t)) \ge 0, \ph t \in [t_1,t_2). \eqno (3.2)
$$
If in addition $\il{t_1}{t}L(s,y_1(s),y_1'(s))ds \ge 0, \ph t \in [t_1,t_2),$ then
$$
\nu(t,y(t),y_1(t),y'(t),y_1'(t)) \ge 0, \ph t \in [t_1,t_2). \eqno (3.3)
$$
}

Proof. Let $\eta(t)$ be a solution of the inequality  (2.11) on $[t_1,t_2)$ with $\eta(t_1) \ge \gamma$ and let $y(t)$ be a solution of Eq. (1.1) with $y(t_1)\in (\gamma, \eta(t_1)], \ph \nu(t_1,\eta(t_1),y(t_1)\eta'(t_1),y'(t_1)) >~ 0, \linebreak \nu(t_1,y(t_1),y_1(t_1),y'(t_1),y_1'(t_1)) \ge 0.$. let $[t_1,t_3)$ be the maximum existence interval for $y(t)$. We must show that
$$
t_3\ge t_2. \eqno (3.4)
$$
Suppose $t_3 < t_2.$ We claim that
$$
y(t) > y_1(t), \ph t \in[t_1,t_3). \eqno (3.5)
$$
Assume this is not true. Then there exists $t_4\in(t_1,t_3)$ such that (since $y(t_1) >  \gamma \ge y_1(t_1)$)
$$
y(t_4) = y_1(t_4). \eqno (3.6)
$$
By (2.3) we have
$$
y(t) - y_1(t) = \exp\biggl\{-\il{t_1}{t}\biggl\{\frac{3}{2} a(\tau)(y(\tau) + y_1(\tau)) + b(\tau)\biggr\} d \tau\biggr\}\times
$$
$$
\times \Biggl[y(t_1) - y_1(t_1) + \il{t_1}{t}\exp\biggl\{\il{t_1}{\tau}\biggl[\frac{3}{2}a(s)(y(s) + y_1(s)) + b(s)\biggr] d s\biggr\} \times
$$
$$
\times\biggl(\nu(t_1,y(t_1),y_1(t_1),y'(t_1),y'_1(t_1)) (\tau - t_1) + \il{t_1}{\tau}[J(s,y(s),y_1(s)) +
$$
$$
+L(s,y_1(s),y_1'(s))]d s\biggr) d \tau\biggr], \phh t \in [t_1,t_3).  \eqno (3.7)
$$
By Lemma 2.3 it follows from the conditions (1) and (2) of the theorem that
$$
y(t) \le \eta(t), \phh t \in [t_1,t_3), \eqno (3.8)
$$
$$
\nu(t,\eta(t),y(t),\eta.(t),y'(t)) \ge 0, \phh t \in[t_1,t_3). \eqno (3.9)
$$
By the mean value theorems for integrals (see [15, p. 869]) it follows from (3.8) that
$$
 \il{t_1}{t}\exp\biggl\{\il{t_1}{\tau}\biggl[\frac{3}{2}a(s)(y(s) + y_1(s)) + b(s)\biggr] d s\biggr\} \times
$$
$$
\times\biggl(\nu(t_1,y(t_1),y_1(t_1),y'(t_1),y'_1(t_1)) (\tau - t_1) + \il{t_1}{\tau}[J(s,y(s),y_1(s)) +
$$
$$
+L(s,y_1(s),y_1'(s))]d s\biggr)d\tau =
$$
$$
= \il{t_1}{\alpha(t)}\exp\biggl\{\il{t_1}{\tau}\biggl[\frac{3}{2}a(s)(\eta(s) + y_1(s)) + b(s)\biggr] d s\biggr\} \times
$$
$$
\times\biggl(\nu(t_1,y(t_1),y_1(t_1),y'(t_1),y'_1(t_1)) (\tau - t_1) + \il{t_1}{\tau}[J(s,y(s),y_1(s)) +
$$
$$
+L(s,y_1(s),y_1'(s))]d s\biggr)d\tau, \phh t\in [t_1,t_4],
$$
for some $\alpha(t) \in [t_1,t], \ph t \in [t_1,t_4]$ By Lemma 2.1 it follows from the conditions (1) and (2) of the theorem that
$$
J(t,y(t),y_1(t)) \ge 0, \phh t \in [t_1,t_4]. \eqno (3.10)
$$
This together with (3.7) and with the condition (3) of the theorem implies that (since $y(t_1) > \gamma) \ph y(t_4)  > y_1(t_4)$, which contradicts (3.5). It follows from (3.5) and from the condition (1) of the theorem that the function $F(t)\equiv \il{t_1}{t}a(\tau) y(\tau) d \tau, \ph t \in [t_1.t_3)$ is bounded from below on $[t_1,t_3)$. In virtue of Lemma 2.2 it follows from here that $[t_1,t_3)$ is not the maximum existence interval for $y(t)$. We have obtained a contradiction, proving (3.4).  From (3.4), (3.5) and (3.9) it follows (3.1) and (3.2). If $\il{t_1}{t}L(s,y_1(s),y_1'(s)ds\ge~ 0, \linebreak t \in[t_1,t_2)$, then by (2,2) from the initial condition $\nu(t_1,y(t_1),y_1(t_1),y'(t_1),y_1'(t_1)) \ge~ 0$ and from (3.10) it follows (3.3). Thus the theorem is proved in the case $y(t_1) > \gamma$. It remains to prove the theorem for the case $y_1(t_1) = \gamma.$ Let for any $\delta > 0$ the function $\widetilde{y}_\delta(t)$ be a solution of Eq. (1.1) with $\widetilde{y}_\delta(t_1) = \gamma +\delta, \ph |\widetilde{y}_\delta'(t_1) - y'(t_1)| <\delta, \ph \nu(t_1,\widetilde{y}_\delta(t_1),y_1(t_1),\widetilde{y}_\delta'(t_1),y'(t_1)) >~ 0$ (the last two inequalities are always satisfiable for all enough small $\delta >0$ due to initial conditions) and let $\eta_\delta(t)$ be a solution of the inequality (2.11) on $[t_1,t_2)$ with $\eta_\delta(t_1) \ge~ \widetilde{y}_\delta(t_1), \linebreak \nu(t_1,\eta_\delta(t_1),\widetilde{y}_\delta(t_1),\eta_\delta'(t_1),\widetilde{y}_\delta'(t_1)) \ge 0$. Then by already proven (by Remark 2.1 if $t_2<\infty$, then $\eta_\delta(t)$ always exists and the proof of the case $t_2 = \infty$ is reducible to the proof of the case $t_2 < \infty$) $\widetilde{y}_\delta(t)$ exists on $[t_1,t_2)$ and
$$
\widetilde{y}_\delta(t) > y_1(t), \phh t \in [t_1,t_2). \eqno (3.11)
$$
Let $[t_1,t_3)$ be the maximum existence interval for $y(t)$. Show that
$$
t_3 \ge t_2. \eqno (3.12)
$$
Suppose $t_3 < t_2$. We claim that
$$
y(t) \ge y_1(t), \phh t \in [t_1,t_3). \eqno (3.13)
$$
Suppose this is not true. Then there exists $t_4\in(t_1,t_3)$ such that
$$
y(t_4) < y_1(t_4). \eqno (3.14)
$$
Since the solutions of Eq. (1.1) are continuously dependent on their ini9tial values we chose $\delta > 0$ so small that
$$
|\widetilde{y}_\delta(t_1) - y(t_1)| < \frac{y_1(t_4) - y(t_4)}{2}.
$$
Then $y(t_4) - y_1(t_4) = y(t_4) - \widetilde{y}_\delta(t_1) + \widetilde{y}_\delta(t_1)  - y_1(t_4) \ge - |y(t_1) - \widetilde{y}_\delta(t_1)| > \frac{y(t_4) - y_1(t_4)}{2}.$ We have obtained a contradiction with (3.14), which proves (3.13). It follows from (3.13) and from the condition (1) of the theorem that the function $F(t)\equiv \il{t_1}{t} a(\tau) d\tau, \ph t \in [t_1,t_3)$ is bounded from below on $[t_1,t_3)$. Then by virtue of Lemma 2.2 $[t_1,t_3)$ is not the maximum existence interval for $y(t)$, which contradicts our supposition. The obtained contradiction proves (3.12). From(3.12) and (3.13) it follows
$$
y(t) \ge y_1(t), \phh t \in [t_1,t_2).
$$
Then by Lemma 2.3 we obtain (3.1) and (3.2). The inequality (3.3) can be proved by analogy with the proof of the already proven  case $y(t_1) > \gamma$. The theorem is proved.

{\bf Remark 3.1.} {\it It is clear from the proof of Theorem 3.1 that the conditions  \linebreak $y(t_1) \in~ [\gamma,\eta(t_1)], \ph \nu(t_1,\eta(t_1),y(t_1),\eta'(t_1),y'(t_1)) >~0$   of Theorem 3.1 can be replaced by the following  ones
$y(t_1) \in (\gamma,\eta(t_1)], \ph \nu(t_1,\eta(t_1),y(t_1),\eta'(t_1),y'(t_1)) \ge~0$.
}

Using Lemma 2.4 instead of Lemma 2.3 by analogy with the proof of Theorem 3.1 one can prove the following theorem

{\bf Theorem 3.2.} {\it Let the following conditions be satisfied

\noindent
 $a(t) > 0, \ph D(t) \ge 0, \ph t \in [t_1,t_2),$

\noindent
 Eq. (2.1) has a solution $y_1(t)$ on $[t_1,t_2)$ such that
$$
\gamma - y_1(t_1) + \il{t_1}{t}\exp\biggl\{\il{t_1}{\tau}\Bigl[\frac{3}{2} a(s)(\zeta(s) + y_1(s)) + b(s)\Bigr]ds\biggr\}\biggl(\il{t_1}{\tau}L(s,y_1(s),y'_1(s)) d s\biggl) d\tau \le 0,
$$
$t \in[t_1,t_2)$ for some $\gamma \le y(t_1)$, where $\zeta(t)$ is a solution of the inequality (2.16) on $[t_1,t_2)$ with $\zeta(t_1) \le \gamma$.

\noindent
Then every solution $y(t)$ of Eq. (1.1) with $y(t_1) \in [\zeta(t_1),\gamma], \ph \nu(t_1,\zeta(t_1),y(t_1),\zeta'(t_1),y'(t_1)) < 0, \ph \nu(t_1,y(t_1),y_1(t_1),y'(t_1),y_1'(t_1)) \le 0$, exists on $[t_1,t_2)$ and
$$
\zeta(t) \le y(t) \le y_1(t), \ph t \in [t_1,t_2),
$$
$$
\nu(t,\zeta(t),y(t),\zeta'(t),y'(t)) \le 0, \ph t \in [t_1,t_2).
$$
If in addition $\il{t_1}{t}L(s,y_1(s),y_1'(s))ds \le 0, \ph t \in [t_1,t_2),$ then
$$
\nu(t,y(t),y_1(t),y'(t),y_1'(t)) \le 0, \ph t \in [t_1,t_2).
$$
}

\phantom{aaaaaaaaaaaaaaaaaaaaaaaaaaaaaaaaaaaaaaaaaaaaaaaaaaaaaaaaaaaaaaaaaaaaaaa} $\blacksquare$

{\bf Remark 3.2.} {\it It is clear from the  Remark 3.1 and the similarity of the proofs of Theorem 3.1 and Theorem 3.2  that the  initial conditions  $y(t_1) \in~ [\zeta(t_1),\gamma], \linebreak \nu(t_1,\zeta(t_1),y(t_1),\zeta'(t_1),y'(t_1)) <~0$   of Theorem 3.2 can be replaced by the following  ones
$y(t_1) \in [\zeta(t_1),\gamma), \ph \nu(t_1,\eta(t_1),y(t_1),\eta'(t_1),y'(t_1)) \le~0$.
}

Let $a_2(t), \ph b_2(t), \ph c_2(t), \ph d_2(t)$ and $e_2(t)$ be real-valued continuous functions on $[t_0,\infty)$. 
We set $L_1(t,u,v)\equiv 3(a_2(t) - a(t))u v + (b_2(t) - b(t))v + (a_2^2(t) - a^2(t)) u^3 + (c_2(t) - c(t)) u^2 + (d_2(t) - d(t)) u + e_2(t) - e(t), \ph t \ge t_0$.
Consider the equation
$$
y'' + 3 a_2(t) y y' + b_2(t) y' + a_2^2(t) y^3 + c_2(t) y^2 + d_2(t) y + e_2(t) = 0, \phh t\in[t_1,2). \eqno (3.15)
$$

{\bf Theorem 3.3.} {\it Let the following conditions be satisfied

\noindent
(4) If $a(t) = 0$, then $c(t) = \frac{3}{2}a'(t), \ph d(t) \ge b'(t)$, otherwise $D(t)\ge 0, \ph t \in[t_1,t_2)$,

\noindent
(5) Eq. (2.1) has a solution $y_1(t)$ on $[t_1,t_2)$ such that $\il{t_1}{t} L(\tau,y_1(\tau),y_1'(\tau)) d\tau \ge 0, \linebreak t \in[t_1,t_2).$

\noindent
(6) Eq. (3.15) has a solution $y_2(t)$ on $[t_1,t_2)$ such that $y_1(t) < y_2(t)$ and \linebreak $\il{t_1}{t} L_1(s,y_2(s),y_2'(s)) d s \ge~ 0, \ph t \in [t_1,t_2).$

\noindent
Then every solution $y(t)$ of Eq. (1.1) with

\noindent
(7) $y_1(t_1) \le y(t_1) \le y_2(t_1),$

\noindent
(8) $\nu(t_1,y(t_1),y_1(t_1),y'(t_1),y_1'(t_1)) \ge 0,$

\noindent
(9) $\nu(t_1,y(t_1),y_2(t_1),y'(t_1),y_2'(t_1)) \le 0,$

\noindent
exists on $[t_1,t_2)$ and
$$
y_1(t) \le y(t) \le y_2(t), \ph [t_1,t_2), \eqno (3.16)
$$
$$
\nu(t,y(t),y_1(t),y'(t),y_1'(t)) \ge 0, \phh  \nu(t,y(t),y_2(t),y'(t),y_2'(t)) \le 0, \phh t \in [t_1,t_2). \eqno (3.17)
$$
}

Proof. Let $y(t)$ be a solution of Eq. (1.1)  with
 $y_1(t_1) < y(t_1) < y_2(t_1)$,
satisfying the initial conditions (8) and (9) and let $[t_1,t_3)$ be its maximum existence interval. We must show
that
$$
t_3 \ge t_2. \eqno (3.18)
$$
Suppose this is not true. We claim that
$$
y_1(t) < y(t) < y_2(t), \ph t \in [t_1,t_3). \eqno (3.19)
$$
Assume at least one of these inequalities is not fulfilled. Then there exists $t_4 \in (t_1,t_3)$ such that
$$
y_1(t) < y(t) < y_2(t), \phh t \in [t_1,t_4), \eqno (3.20)
$$
$$
y_1(t_4) = y(t_4) \ph \mbox{or} \ph y(t_4) = y_2(t_4). \eqno (3.21)
$$
By (2.3) we have
$$
y(t_4) - y_1(t_4) = \exp\biggl\{-\il{t_1}{t_4}\biggl\{\frac{3}{2} a(\tau)(y(\tau) + y_1(\tau)) + b(\tau)\biggr\} d \tau\biggr\}\times
$$
$$
\times \Biggl[y(t_1) - y_1(t_1) + \il{t_1}{t_4}\exp\biggl\{\il{t_1}{\tau}\biggl[\frac{3}{2}a(s)(y(s) + y_1(s)) + b(s)\biggr] d s\biggr\} \times
$$
$$
\times\biggl(\nu(t_1,y(t_1),y_1(t_1),y'(t_1),y'_1(t_1)) (\tau - t_1) + \il{t_1}{\tau}[J(s,y(s),y_1(s)) +
$$
$$
+L(s,y_1(s),y_1'(s))]d s\biggr) d \tau\biggr],   \eqno (3.22)
$$
$$
y_2(t_4) - y(t_4) = \exp\biggl\{-\il{t_1}{t_4}\biggl\{\frac{3}{2} a(\tau)(y_2(\tau) + y(\tau)) + b(\tau)\biggr\} d \tau\biggr\}\times
$$
$$
\times \Biggl[y_2(t_1) - y(t_1) + \il{t_1}{t_4}\exp\biggl\{\il{t_1}{\tau}\biggl[\frac{3}{2}a(s)(y_2(s) + y(s)) + b(s)\biggr] d s\biggr\} \times
$$
$$
\times\biggl(\nu(t_1,y_2(t_1),y(t_1),y_2'(t_1),y'(t_1)) (\tau - t_1) + \il{t_1}{\tau}[J(s,y_2(s),y(s)) -
$$
$$
-L_1(s,y_2(s),y_2'(s))]d s\biggr) d \tau\biggr],   \eqno (3.23)
$$
It follows from the condition (4) of the theorem  and (3.20) that
$$
J(t,y(t),y_1(t)) \ge 0, \phh J(t,y_2(t),y(t)) \le 0, \phh t \in [t_1,t_4].
$$
These inequalities together with (3.22), (3.23) and with conditions (5), (6) of the theorem imply that $y_1(t_4) < y(t_4) < y_2(t_4)$, which contradicts (3.21). The obtained contradiction proves (3.19). It follows from (3.19) that the function $ y(t),  \ph t \in [t_1,t_3)$ is bounded on $[t_1,t_3)$. By virtue of Lemma 2.5 it follows from here that $[t_1,t_3)$ is not the maximum existence interval for $y(t)$, which contradicts our supposition. The obtained contradiction proves (3.18). From (3.18) it follows (3.16). By (2.2) we have
$$
\nu(t,y(t),y_1(t),y'(t),y_1'(t)) = \nu(t_1,y(t_1),y_1(t_1),y'(t_1),y_1'(t_1))  +
$$
$$
+\il{t_1}{t}[J(\tau,y(\tau),y_1(\tau)) + L(\tau,y_1(\tau),y_1'(\tau))]d\tau, \phh t \in [t_1,t_2), \eqno (3.24)
$$
$$
\nu(t,y_2(t),y(t),y_2'(t),y'(t)) = \nu(t_1,y_2(t_1),y(t_1),y_2'(t_1),y'(t_1))  +
$$
$$
+\il{t_1}{t}[J(\tau,y_2(\tau),y(\tau)) + L_1(\tau,y_2(\tau),y_2'(\tau))]d\tau, \phh t \in [t_1,t_2), \eqno (3.25)
$$
It follows from (4) and (3.16) that
$$
J(t,y(t),y_1(t)) \ge 0, \phh J(t,y_2(t),y(t)) \le 0, \phh t\in [t_1,t_2).
$$
These inequalities with (3.14), (3.25) and with the conditions (5), (6) imply (3.17). Thus the theorem is proved in the particular case of $y_1(t_1) < y(t_1) < y_2(t_1)$. Let $y(t_1) = y(t_1)$ and let $[t_1,t_3)$ be the maximum existence interval for $y(t)$. Show that
$$
t_3 \ge t_1. \eqno (3.26)
$$
Suppose $t_3 < t_2$. We claim that
$$
y_1(t) \le y(t) \le y_2(t), \phh t \in [t_1,t_3). \eqno (3,27)
$$
Let $\widetilde{y}_\delta(t)$ be a solution of Eq. (1.1) with $\widetilde{y}_\delta(t_1) = y_1(t_1) + \delta, \ph |\widetilde{y}_\delta '(t_1) - y'(t_1)| < \delta,$ where $\delta > 0$ is so small  that $y_1(t_1) +\delta < y_2(t_1), \ph \nu(t,\widetilde{y}_\delta(t_1),y_1(t_1),\widetilde{y}_\delta'(t_1),y_1'(t_1)) \ge 0, \ph \nu(t,\widetilde{y}_\delta(t_1),y_2(t_1),\widetilde{y}_\delta'(t_1),y_2'(t_1)) \le 0.$ By already proven
$$
y_1(t) < \widetilde{y}_\delta(t) < y_2(t), \ph t \in [t_1,t_2). \eqno (3.28)
$$
Using (2.3) one can show that
$$
y(t) < y_2(t), \phh t \in [t_1,t_3). \eqno (3.29)
$$
Suppose (3.27) is not true. Then taking into account (3.29) we conclude that
$$
y(t_4) < y_1(t_4) \eqno (3.30)
$$
for some $t_4 \in (t_1,t_3)$. Since the solutions of Eq. (1.1) are continuously dependent on their initial values we chose $\delta > 0$ so small that $|\widetilde{y}_\delta(t_4) - y(t_4)| < \frac{y_1(t_4) - y(t_4)}{2}.$ Then $|y(t_4) - y_1(t_4)| = y(t_4) - \widetilde{y}_\delta(t_4) + \widetilde{y}_\delta(t_4) - y_1(t_4) \ge -|y(t_4) - \widetilde{y}_\delta(t_4)| > \frac{y(t_4) - y_1(t_4)}{2}.$ We have obtained a contradiction with (3.30). The obtained contradiction proves (3.27). By virtue of Lemma 2.5 it follows  from (3.27) that $[t_1,t_3)$ is not the maximum existence interval for $y(t)$, which contradicts our supposition. The obtained contradiction proves (3.26). From (3.26) and (3.27) it follows (3.16). The proof of (3.16) in the case $y(t_1) = y_2(t_1)$ can be realized by analogy with the proof of the case $y(t_1) = y_1(t_1)$. Then using (2.2) on the basis of (3.16) one can prove (3.17). The theorem is proved.

\vskip 10pt

{\bf 4. Global solvability criteria.}

{\bf Theorem 4.1.} {\it Let the following conditions be satisfied

\noindent
$(\alpha) \ph a(t) > 0, \ph D(t) \ge 0, \ph t \ge t_0,$

\noindent
$(\beta) \ph M\equiv \sup\limits_{t \ge t_0} \frac{|c(t)| + |d(t)| + |e(t)|}{a^2(t)} < +\infty,$

\noindent
$(\gamma) \ph \il{t_0}{t}\exp\biggl\{\il{t_0}{\tau}\Bigr[\frac{3}{2} M a(s) + b(s)\Bigr]d\tau\biggr\}\Biggl(\il{t_0}{\tau} e(s) d s\Biggr) d \tau \le 0, \ph t \ge t_0.$

\noindent
Then every solution $y(t)$ of Eq. (1.1) with $y(t_0) \in [0,M], \ph \nu(t_0,M,y(t_0),0,y'(t_0)) \ge~ 0, \linebreak \nu(t_0,0,y(t_0),0,y'(t_0)) \le 0$ exists on $[t_0,\infty)$ and
$$
0\le y(t) \le M, \phh t \ge t_0, \eqno (4.1)
$$
$$
\nu(t,M,y(t),0,y'(t)) \ge 0, \ph t \ge t_0. \eqno (4.2)
$$
$$
\nu(t,0,y(t),0,y'(t)) \le 0,  \ph t \ge t_0. \eqno (4.3)
$$
}

Proof. By Remark 2.1 it follows from the conditions $(\alpha)$ and $(\beta)$ that $\eta(t) \equiv M, \ph t \ge t_0$ is a solutions of the inequality (2.11) on $[t_0,\infty)$. We put $a_1(t) = a(t), \ph b_1(t) = b(t), \ph c_1(t) = c(t), \ph d_1(t) = d(t), \ph e_1(t) \equiv 0, \ph t \ge t_0$. in Eq. (2.1). Then, obviously, $y_1(t)\equiv 0$ is a solution of Eq. (2.1) on $[t_0,\infty)$. Then by Theorem 3.1
$$
\il{t_0}{t}\exp\biggl\{\il{t_0}{\tau}\Bigl[\frac{3}{2} M a(s) + b(s)\Bigr]d s\biggr\}\biggl(\il{t_0}{\tau}L(s,0,0)d s\biggr) d \tau \ge 0, \phh t \ge t_0 \eqno (4.11)
$$
and the condition $(\alpha)$ holds, then every solution $y(t)$ of Eq. (1.1) with $y(t_0) \in [0,M], \linebreak \nu(t,M,y(t),0,y'(t)) \ge 0, \ph \nu(t,0,y(t),0,y'(t)) \le 0$ exists on $[t_0,\infty)$ and the inequalities (4.1)-(4.3) hold. Since $L(t,0,0) = - e(t), \ph t \ge t_0$ it follows from the condition $(\gamma)$ the inequality (4.1). The theorem is proved.

Using Theorem 3.2 instead of Theorem 3.1 by analogy with the proof of Theorem 3.1 one can prove the following theorem

{\bf Theorem 4.2.} {\it Let the following conditions be satisfied

\noindent
$ a(t) > 0, \ph D(t) \ge 0, \ph t \ge t_0, \ph M\equiv \sup\limits_{t \ge t_0} \frac{|c(t)| + |d(t)| + |e(t)|}{a^2(t)} < +\infty,$

\noindent
$ \il{t_0}{t}\exp\biggl\{\il{t_0}{\tau}\Bigr[\frac{3}{2} M a(s) + b(s)\Bigr]d\tau\biggr\}\Biggl(\il{t_0}{\tau} e(s) d s\Biggr) d \tau \le 0, \ph t \ge t_0.$

\noindent
Then every solution $y(t)$ of Eq. (1.1) with $y(t_0) \in [-M,0], \ph \nu(t_0,M,y(t_0),0,y'(t_0)) \le~ 0, \linebreak \nu(t_0,0,y(t_0),0,y'(t_0)) \ge 0$ exists on $[t_0,\infty)$ and
$$
- M\le y(t) \le 0, \ph
\nu(t,M,y(t),0,y'(t)) \le 0, \ph
\nu(t,0,y(t),0,y'(t)) \ge 0,  \ph t \ge t_0.
$$
}

\phantom{aaaaaaaaaaaaaaaaaaaaaaaaaaaaaaaaaaaaaaaaaaaaaaaaaaaaaaaaaaaaaaaaaaaaaaa} $\blacksquare$

Let $t_0 < t_1 \ldots < t_n < \ldots$ be an infinitely large sequence. We set
$$
M_n\equiv \max_{t\in[t_0,t_n]}\frac{|c(t)| + |d(t)| = |e(t)|}{a^2(t)|}, \phh n=1,2, \ldots.
$$
{\bf Theorem 4.3} {\it Let the following conditions be satisfied.

\noindent
$(\alpha) \ph a(t) > 0, \ph D(t) \ge 0, \ph t \ge t_0,$

\noindent
$(\delta) \ph \il{t_n}{t}\exp\biggl\{\il{t_n}{\tau}\Bigl[\frac{3}{2}a(s)M_{n+1} + b(s)\Bigr]ds\biggr\}\biggl(\il{t_n}{\tau}e(s) d s\biggr) d\tau \le 0, \ph t \in [t_n,t_{n+1}), \ph n=0,1,2, \ldots.$

\noindent
Then every solution $y(t)$ of Eq. (1.1) with $y(t_0) \in [0,M_1], \ph \nu(t_0,y(t_0),M_1,y'(t_0),0) \le 0, \ph \nu(t_0,y(t_0),0,y'(t_0),0) \ge 0$ exists on $[t_0,\infty)$ and
$$
0\le y(t) \le M_n, \phh t \in [t_{n-1},t_n], \eqno (4.5)
$$
$$
\nu(t,y(t),M_n,y'(t),0) \le 0, \ph \nu(t,M_n,y(t),0,y'(t)) \ge 0, \ph t\in [t_{n-1},t_n], \ph n=1,2,\ldots. \eqno (4.6)
$$
}

Proof. According to Remark 2.1 $\eta_n(t)\equiv M_n, \ph t \in[t_{n-1},t_n]$ is a solution of the inequality (2.11) on $[t_{n-1},t_n], \ph n=1,2,\ldots.$ We put $a_1(t)=a(t), \ph b_1(t) = b(t), \ph c_1(t)=c(t),\ph d_1(t)=d(t), \ph e_1(t)\equiv0, \ph t \ge t_0$ in Eq. (2.1). In this case, obviously, $y_1(t)\equiv 0$ is a solution of Eq. (2.1) on $[t_0,\infty)$. Then by Theorem 3.1 if
$$
\il{t_n}{t}\exp\biggl\{\il{t_n}{\tau}\Bigl[\frac{3}{2} a(s)M_{n+1} = b(s)\Bigr]d s\biggr\}\biggl(\il{t_n}{\tau}L(s,0,0)ds\Biggr)d\tau \ge 0, \ph t\in [t_n,t_{n+1}],  \eqno (4.7)
$$
$n=0,1,2,\ldots$ every solution $y_n(t)$ of Eq. (1.1) with $y_n(t_{n-1}) \in [0,M_n], \linebreak \nu(t_{n-1},y_n(t_{n-1}),0,y_n'(t_{n-1}),0) \ge 0, \phh \nu(t_{n-1},M_n,y_n(t_{n-1}),0,y_n'(t_{n-1}),) \ge 0$ exists on $[t_{n-1},t_n]$ and
$$
0 \le y_n(t) \le M_n, \phh t \in [t_{n-1},t_n],  \ph n=1,2,\ldots. \eqno (4.8)
$$
$$
\nu(t,y_n(t),0,y_n'(t),0) \ge 0, \ph \nu(t,M_n,y_n(t),0) \ge 0, \ph t \in [t_{n-1},t_n], \ph n=1,2,\ldots.\eqno (4.9)
$$
Since $M_1 \le M_2 \le \ldots$ according to (4.8) we can take that $y_{n+1}(t_n) = y_n(t_n), \ph y_{n+1}'(t_n) = y_n'(t_n), \ph n=1,2,\ldots.$ Note that according to (4.9) in this particular case of chose of initial conditions the relations ( initial conditions) $\nu(t_{n-1},M_n,y_{n-1}(t_{n-1}),0,y_n'(t_{n-1})) \ge 0, \linebreak \nu(t_{n-1},y_{n-1}(t_{n-1}),0,y_n'(t_{n-1}),0) \ge 0, \ph n=2,3\ldots$ remain valid. Hence the function
$$
y(t) = y_n(t), t\in [t_{n-1},t_n], \phh n=1,2,\ldots
$$
is a solution of Eq. (1.1) on $[t_0,\infty)$, satisfying the initial conditions $y(t_0) \in [0,M_1], \linebreak \nu(t_{0},M_n,y(t_0),0,y'(t_{n-1})) \ge 0, \ph \nu(t_0,y(t_0),M,y'(t_0),0) \le 0$, for which according to (4.8) and (4.9) the inequalities hold. Then since $L(t,y_1(t),y_1'(t)) = L(t,0,0) = - e(t), \ph t \ge t_0$ the inequality (4.7) follows from the condition $(\delta)$. The theorem is proved.

Using Theorem 3.2 instead of Theorem 3.1 by analogy with the proof of Theorem 4.3  can be proved the following theorem.

{\bf Theorem 4.4.} {\it Let the following conditions be satisfied

\noindent
$a(t) < 0, \ph D(t) \ge 0, \ph t \ge t_0,$

\noindent
$\il{t_n}{t}\exp\biggl\{-\il{t_n}{\tau}\Bigl[\frac{3}{2}a(s) M_{n+1} + b(s)\Bigr]d s\biggr\}\biggl(\il{t_n}{\tau}e(s) d s\biggr) d\tau \ge 0, \ph t \in [t_n,t_{n-1}), \ph n=0,1,2,\ldots.$

\noindent
Then every solution $y(t)$ of Eq. (1.1) with $y(t) \in [-M_1,0], \ph \nu(t_0,y(t_0),0,y'(t_0)) \le 0, \linebreak \nu(t_0,-M_1,y(t_0),0,y'(t_0)) \le 0$ exists on $[t_0,\infty)$ and
$$
-M_n \le y(t) \le 0, \phh \nu(t,y(t),0,y'(t),0) \le 0, \phh \nu(t,-M_n,y(t),0y'(t)) \le 0, \ph t \in [t_{n-1},t_n],
$$
$n=1,2,\ldots$.
}

\phantom{aaaaaaaaaaaaaaaaaaaaaaaaaaaaaaaaaaaaaaaaaaaaaaaaaaaaaaaaaaaaaaaaaaaaaaa} $\blacksquare$

We set 
$$
\rho_\pm(t) \equiv \frac{- d(t) \pm\sqrt{d^2(t) - 4 c(t) e(t)}}{2 c(t)}, \phh t \ge t_0. 
$$

{\bf Theorem 4.5.} {\it Let the conditions (4) of Theorem 3.3 and the conditions

\noindent
(10) $c(t) \ne 0, \ph d^2(t) - 4 c(t) e(t) > 0, \ph \rho_\pm(t) \in C^2([t_0,\infty)), \ph  t \ge t_0,$

\noindent 
(11) $\il{t_0}{t}[\rho_-''(\tau) + 3 a(\tau) \rho_-(\tau)\rho_-'(\tau) + a^2(\tau)\rho_-^3(\tau)] d \tau \le 0, \ph t \ge t_0,$

\noindent
(12) $\il{t_0}{t}[\rho_+''(\tau) + 3 a(\tau) \rho_+(\tau)\rho_+'(\tau) + a^2(\tau)\rho_+^3(\tau)] d \tau \ge 0, \ph t \ge t_0$

\noindent
be satisfied.

\noindent
Then every solution $y(t)$ of Eq. (1.1) with the initial  conditions $y(t_0)\in [\rho_-(t_0),\rho_+(t_0)], \linebreak \nu(t_0,y(t_0),\rho_-(t_0),y'(t_0),\rho_-'(t_0)) \ge 0, \ph \nu(t_0,y(t_0),\rho_+(t_0),y'(t_0),\rho_+'(t_0)) \le 0$ exists on $[t_0,\infty)$ and
$$
\rho_-(t) \le y(t) \le \rho_+(t) ,\ph t \ge t_0. \eqno (4.10)
$$
$$
\nu(t,y(t),\rho_-(t),y'(t),\rho_-'(t)) \ge 0, \phh \nu(t,y(t),\rho_+(t),y'(t),\rho_+'(t)) \le 0,. \ph t \ge t_0. \eqno (4.11)
$$
}

Proof. We put $a_1(t) = a_2(t) = a(t), \ph b_1(t) = b_2(t) = b(t), \ph c_1(t) = c_2(t) = c(t), \ph d_1(t) = d_2(t) = d(t), \ph e_1(t) = e(t)  -\rho_-''(t) - 3 a(\tau) \rho_-(t)\rho_-'(t) - a^2(t)\rho_-^3(t), \ph e_2(t) =e(t)  -\rho_-''(t) - 3 a(\tau) \rho_-(t)\rho_-'(t) - a^2(t)\rho_-^3(t), \ph t \ge t_0.$ Then it is not difficult to verify that $\rho_-(t)$ is a solution of Eq. (2.1) on $[t_0,\infty)$ and $\rho_+(t)$ is a solution of Eq. (3.15) on $[t_0,\infty)$. Then the conditions (4), (10)-(12) provide the satisfiability all of the conditions of Theorem 3.3. Hence, every solution $y(t)$ of Eq. (1.1)  with the initial conditions  $y(t_0)\in [\rho_-(t_0),\rho_+(t_0)], \linebreak \nu(t_0,y(t_0),\rho_-(t_0),y'(t_0),\rho_-'(t_0)) \ge 0, \ph \nu(t_0,y(t_0),\rho_+(t_0),y'(t_0),\rho_+'(t_0)) \le 0$ exists on $[t_0,\infty)$  and the relations (4.10) and (4.11) are fulfilled. The theorem is proved.

\vskip 10pt

{\bf 5. An application to three dimensional linear systems of ordinary differential equations.} Let $a_{jk}(t), \ph j,k=\overline{1,3}$ be real-valued continuous functions on $[t_0,\infty)$. Consider the linear system of ordinary differential equations
$$
\sisttt{\phi' = a_{11}(t)\phi + a_{12}(t)\psi + a_{13}(t) \chi,}{\psi' = a_{21}(t)\phi + a_{22}(t)\psi + a_{23}(t) \chi,}{\chi' = a_{31}(t)\phi + a_{32}(t)\psi + a_{33}(t) \chi, \ph t \ge t_0.} \eqno (5.1)
$$
Assume
$$
a_{13}(t) \equiv 0, \phh t \ge t_0. \eqno (5.2)
$$
Let us substitute
$$
\psi = y \phi \eqno (5.3)
$$
in the system (5.1). We obtain
$$
\sisttt{\phi' = [a_{11}(t) + a_{12}(t)y],}{y\phi' = [a_{21}(t) - y' + a_{22}(t) y]\phi + a_{23}(t)\chi,}{\chi' = [a_{31}(t) + a_{32}(t) y] \phi + a_{33}(t) \chi, \ph t \ge t_0.} \eqno (5.4)
$$
From the first and second equations of the last system we obtain
$$
[y' + a_{12}(t) y^2 + B(t) y - a_{21}(t)]\phi = a_{23}(t) \chi, \phh t \ge t_0, \eqno (5.5)
$$
where $B(t) \equiv a_{11}(t) - a_{22}(t), \ph t \ge t_0.$ Assume
$$
a_{23}(t) \ne 0 ,\phh t \ge t_0. \eqno (5.6)
$$
Then from (5.5) we get
$$
\chi = \frac{y' + a_{12}(t) y^2 + B(t) y - a_{21}(t)}{a_{23}(t)} \phi, \phh t \ge t_0. \eqno (5.7)
$$
This together with the first and the third equations of the system (5.4) implies
$$
\biggl\{\frac{y' + a_{12}(t)y^2 + B(t) y - a_{21}(t)}{a_{23}(t)}\biggr\}' \phi +
$$
$$
+\biggl\{\frac{y' + a_{12}(t)y^2 + B(t) y - a_{21}(t)}{a_{23}(t)}\biggr\}[a_{11}(t) + a_{12}(t) y]\phi =
$$
$$
=\biggl[a_{31}(t) + a_{32}(t) y + a_{33}(t)\frac{y' + a_{12}(t) y^2 + B(t) y - a_{21}(t)}{a_{23}(t)} \biggr] \phi, \ph t \ge t_0.
$$
By (5.3) and (5.7) it follows from here that under the restrictions (5.2), (5.6) all solutions  $y(t)$ of the equation
$$
\biggl\{\frac{y' + a_{12}(t)y^2 + B(t) y - a_{21}(t)}{a_{23}(t)}\biggr\}' +
$$
$$
 +[a_{11}(t) - a_{33}(t) + a_{12}(t) y] \biggl\{\frac{y' + a_{12}(t)y^2 + B(t) y - a_{21}(t)}{a_{23}(t)}\biggr\}- a_{32}(t) y - a_{31}(t) = 0,  \eqno (5.8)
$$
$t \ge t_0,$ existing on any interval $[t_1.t_2) \subset [t_0,\infty)$, are connected with solutions $(\phi(t),\psi(t),\chi(t))$ of the system (5.1) by relations
$$
\phi(t) = \phi(t_1)\exp\biggl\{\il{t_1}{t}[a_{11}(\tau) + a_{12}(\tau) y(\tau)]d\tau\biggr\}, \ph  \phi(t_1) \ne 0, \ph t\in [t_1.t_2), \eqno (5.9)
$$
$$
\psi(t) = y(t) \phi(t), \ph \chi(t) = \frac{y'(t) + a_{12}(t) y^2(t) + B(t) y(t) - a_{21}(t)}{a_{23}(t)} \phi(t), \ph t\in [t_1.t_2). \eqno (5.10)
$$
Assume
$$
B(t), \ph a_{12}(t), \ph a_{23}(t) \in C^1([t_0,\infty)). \eqno (5.11)
$$
Then making some differentiations and arithmetic simplifications in (5.8) we obtain
$$
y'' + 3 a_{12}(t) y y' + A(t) y' + a_{12}^2(t) y^3 + C(t) y^2 + F(t) y + E(t) = 0, \phh t \ge t_0, \eqno (5.12)
$$
where
$$
A(t) \equiv 2 a_{11}(t) - a_{22}(t) - a_{33}(t) -  \frac{a_{23}'(t)}{a_{23}(t)},
$$
$$
C(t) \equiv a_{12}(t)[2a_{11}(t) - a_{22}(t) - a_{33}(t)] + a_{23}(t)\biggl(\frac{a_{12}(t)}{a_{23}(t)}\biggr)',
$$
$$
F(t) \equiv [a_{11}(t) - a_{22}(t)] [a_{11}(t) - a_{33}(t)] - a_{23}(t) a_{32}(t) + a_{23}(t)\biggl(\frac{a_{11}(t) - a_{22}(t)}{a_{23}(t)}\biggr)',
$$
$$
E(t)\equiv[a_{33}(t) - a_{11}(t)] a_{21}(t) - a_{23}(t) a_{31}(t) - a_{23}(t)\biggl(\frac{a_{21}(t)}{a_{23}(t)}\biggr)', \phh  t \ge t_0.
$$

{\bf Remark 5.1.} {\it If $\frac{a_{13}(t)}{a_{12}(t)}$ is well defined on $[t_0,\infty)$ and $\in C^1([t_0,\infty))$ (that is $a_{13}(t) = a_{12}(t)\lambda(t), \ph \lambda(t) \in C^1([t_0,\infty))$), then the linear transformation
$$
\psi = \eta - \lambda(t)\chi
$$
reduces the general system (5.1) to its particular case of (5.2) (of $a_{13}(t)\equiv 0, \ph t\ge t_0.$)
}

{\bf Definition 5.1.} {\it The system (5.1) is called non oscillatory, if it has a solution \linebreak $(\phi(t),\psi(t),\chi(t))$ such that $\phi(t) \ne 0, \ph t \ge t_0.$
}

By (5.9) from the results of the section 4 we obtain

{\bf Theorem 5.1.} {\it Let the conditions (5.2), (5.6) and (5.11) be satisfied and let for $a(t) =~ a_{12}(t), \ph b(t) = A(t), \ph c(t) = C(t), \ph d(t) = D(t), \ph e(t) = E(t), \ph t \ge t_0$ the conditions of one of Theorems 4.1 - 4.5 be satisfied. Then the system (5.1) is non oscillatory.
}

\vskip 20pt

\centerline{\bf References}

\vskip 20pt

\noindent
1.  G. A. Grigorian,  On two comparison tests for second-order linear ordinary differential \linebreak \phantom{a}
equations, Diff. Urav., vol 47, (2011), 1225--1240 (in Russian), Diff. Eq., vol. 47, (2011),   \linebreak \phantom{a}
1237--1252  (in English).

\noindent
2.  G. A. Grigorian,   Two comparison criteria for scalar Riccati equations and some of  \linebreak \phantom{a} their applications
Izv. Vyssh. Uchebn. Zaved. Mat.(2012), no. 11, 20–35.
Russian Math. (Iz. VUZ) 56 (2012), no. 11, 17–30.

\noindent
3.   G. A. Grigorian,   Interval oscillation criteria for linear matrix Hamiltonian systems.  \linebreak \phantom{a}
Rocky Mountain J. Math. 50 (2020), no. 6, 2047–2057.

\noindent
4.  G. A. Grigorian, New reducibility criteria for systems of two linear first-order ordinary  \linebreak \phantom{a} differential equations.
Monatsh. Math. 198 (2022), no. 2, 311–322.

\noindent
5.  G. A. Grigorian, Oscillation and non-oscillation criteria for linear nonhomogeneous  \linebreak \phantom{a} systems of two first-order ordinary differential equations.
J. Math. Anal. Appl. 507  \linebreak \phantom{a} (2022), no. 1, Paper No. 125734, 10 pp.

\noindent
6.  G. A. Grigorian, Oscillatory and non-oscillatory criteria for linear four-dimensional   \linebreak \phantom{a} Hamiltonian systems
Math. Bohem. 146 (2021), no. 3, 289–304.

\noindent
7.  G. A. Grigorian,  On the reducibility of systems of two linear first-order ordinary  \linebreak \phantom{a} differential equations.
Monatsh. Math. 195 (2021), no. 1, 107–117.

\noindent
8.  G. A. Grigorian, Oscillation criteria for linear matrix Hamiltonian systems.
Proc. Amer.  \linebreak \phantom{a} Math. Soc. 148 (2020), no. 8, 3407–3415.

\noindent
9.  G. A. Grigorian,  Oscillatory and non oscillatory criteria for the systems of two linear  \linebreak \phantom{a} first order two by two dimensional matrix ordinary differential equations.
Arch. Math. \linebreak \phantom{a} (Brno) 54 (2018), no. 4, 189–203.

\noindent
10.  G. A. Grigorian,   Stability criterion for systems of two first-order linear ordinary   \linebreak \phantom{a}  differential equations.
Mat. Zametki 103   (2018), no. 6, 831–840.
Math. Notes 103  \linebreak \phantom{a} (2018), no. 5-6, 892–900.

\noindent
11. G. A. Grigorian,  Oscillatory criteria for the systems of two first-order linear differential  \linebreak \phantom{a} equations.
Rocky Mountain J. Math. 47 (2017), no. 5, 1497–1524.

\noindent
12.  G. A. Grigorian, Some properties of the solutions of third order linear ordinary   \linebreak \phantom{a} differential equations
Rocky Mountain J. Math. 46 (2016), no. 1, 147–168.

\noindent
13.  G. A. Grigorian, Stability criteria for systems of two first-order linear ordinary differential  \linebreak \phantom{a}  equations.
Math. Slovaca 72 (2022), no. 4, 935–944.

\noindent
14. B. P. Demidovich, Lectures on the mathematical theory of stability. Moskow, \linebreak \phantom{a} "Nauka", 1967.

\noindent
15. R. E. Edvards, A Formal Background to Mathematics (Springer-Verlag, New York,  \linebreak \phantom{a}  Heidelberg, Berlin, 1980).  MR0636095

\end{document}